\documentclass[11pt]{article}
\usepackage{amssymb}
\usepackage{mathtools}
\usepackage{mathrsfs}
\usepackage{amsmath}
\usepackage{amsfonts}
\usepackage{stmaryrd}
\usepackage{amsthm}
\usepackage{thmtools}
\usepackage{thm-restate}

\usepackage{thmtools}
\usepackage{tikz-cd}

\usepackage{enumitem,kantlipsum}
\input amssym.def
\input amssym
\input xypic
\usepackage{comment}
\def\demo{\noindent{\bf Proof. }}
\def\QED{\hfill$\Box$}

\newtheorem*{Theorem*}{Theorem}
\newtheorem{Theorem}{Theorem}[section]
\newtheorem{Lemma}[Theorem]{Lemma}

\newtheorem{Proposition}[Theorem]{Proposition}
\newtheorem{Example}[Theorem]{Example}

\newtheorem{rem}[Theorem]{Remark}
\newtheorem{conv}[Theorem]{Convention}

\begin{document}
\title{
{{\bf On the existence of birational maximal Cohen-Macaulay modules over biradical extensions in mixed characteristic}} \\
}

\author{ Prashanth Sridhar\footnote{Corresponding Author}\\
\begin{tabular}{c}
\vspace{0.15in} \      \\
University of Kansas\\ 405 Snow Hall, 1460 Jayhawk Blvd,  \\
Lawrence, Kansas, USA 66045 \\
e-mail: prashanth@ku.edu
\end{tabular}}

\date{ }

\maketitle \vspace{-0.2in}

\begin{abstract}
\noindent Let $S$ be an unramified regular local ring of mixed characteristic $p\geq 3$ and $S^p$ the subring of $S$ obtained by lifting to $S$ the image of the Frobenius map on $S/pS$. Let $R$ be the integral closure 
of $S$ in a biradical extension of degree $p^2$ of its quotient field obtained by adjoining $p$-th roots of sufficiently general square free elements $f,g\in S^p$. 
We show that $R$ admits a birational maximal Cohen-Macaulay module. It is noted that $R$ is not automatically Cohen-Macaulay. 
\\
\textbf{Classification Codes}: 13B22, 13C10, 13C14, 13C15, 13D22, 13E05, 13H05.
\end{abstract}

\section{Introduction}
The existence of maximal Cohen-Macaulay (MCM) modules (\ref{notation}) over catenary rings is largely an open problem. Hochster conjectured that every complete local domain admits an MCM module (see for example \cite{10.1007/BFb0068925} and \cite{R}), but this is known to be true only in very few cases. The primary goal of this article is to show the existence of an MCM module for a new class of rings.
\par Hochster's conjecture reduces to the integral closure of a complete regular local ring in a normal extension of its fraction field, so that one may approach the problem from this viewpoint. Let $S$ denote an unramified regular local ring with quotient field $L$ and $R$ the integral closure of $S$ in a finite field extension $K/L$. In \cite{RO} it is shown that if $K/L$ is Abelian and the degree of the extension is not divisible by the characteristic of the residue field, then $R$ is Cohen-Macaulay. Note that this applies when $S$ contains the rational numbers. The proof uses the fact that the group algebra $k[G]$ is a product of fields when the residue field $k$ of $S$ is algebraically closed and $G=Gal(K/L)$. There is no direct analog of the argument when $char(k)$ divides the order of $G$. In fact, the conclusion fails in mixed characteristic as shown in \cite{KO} and \cite{KA}.

\par Motivated by the above phenomenon and the studies in \cite{KA}, \cite{KATZ2021350} and \cite{PS}, we consider extensions $K/L$ with the property that $R$ is $S$-free when $S$ contains a field but not necessarily so otherwise. Kummer theory tells us that Abelian extensions of a field of characteristic zero containing ``suitable" roots of unity are repeated radical extensions. Thus it is natural to study repeated $n$-th root extensions of an unramified regular local ring of mixed characteristic $p>0$ with the property that $p|n$. If $S$ were to contain the rational numbers or if $S$ were of mixed characteristic $p>0$ and $p\nmid n$, then it follows that the integral closure of $S$ in an arbitrary repeated radical extension is Cohen-Macaulay, see \cite{PS} and \cite{HK}. 
\par Now assume additionally that $S$ has mixed characteristic $p\geq 3$. If $K$ is the extension by a $p$-th root of a square free element of $S$, then $R$ is Cohen-Macaulay as shown in \cite{KA}. In contrast, we shall see that the integral closure $R$ fails to be Cohen-Macaulay in a finite square free tower of $p$-th roots. In this paper, we consider the existence of maximal Cohen-Macaulay modules over the integral closure of $S$ in an extension of degree $p^2$, obtained by adjoining $p$-th roots of sufficiently general square free elements. The generality and square free conditions on the elements are natural, since any given multi-radical extension can be embedded in a sufficiently large, general, square free tower.
\par Let $S^p$ denote the subring of $S$ obtained by lifting the image of the Frobenius map on $S/pS$ to $S$. Towards the above goal, the immediate obstruction is when $f,g\in S^p$, see \cite{DKPS}. We now outline our principal findings. Fix $f,g \in S^p$ square free, non-units that form a regular sequence in $S$ or units that are not $p$-th powers in $S$. Let $\omega^p=f$ and $\mu^p=g$. If $f,g$ are units, assume further that $[L(\omega,\mu):L]=p^2$. Given integers $n,k\geq 1$, let $S^{p^k\wedge p^n}\subset S$ be the multiplicative subset of $S$ consisting of elements expressible in the form $x^{p^k}+y\cdot p^n$ for some $x,y\in S$. For a discussion of the case $p=2$, see \cite{PS}, where the results are sharper since in that case such extensions are automatically Abelian. When $p\geq 3$, the presence of a $p$-th root of unity in $S$ necessarily ramifies $p$, so we do not quite have the same leg room. The main result of this paper is
\begin{Theorem*}\noindent$\mathbf{4.10}$
Let $(S,\mathfrak{m})$ be an unramified regular local ring of mixed characteristic $p\geq 3$. Then
\begin{enumerate}
    \item $R$ is Cohen-Macaulay if
    \begin{enumerate}
        \item At least one of $S[\omega], S[\mu]$ is not integrally closed.
        \item $S[\omega],S[\mu]$ are integrally closed and $fg^i\notin S^{p\wedge p^2}$ for all $1\leq i\leq p-1$.
    \end{enumerate}
    \item Let $S[\omega],S[\mu]$ be integrally closed such that $fg\in S^{p\wedge p^2}$. Then $R$ is Cohen-Macaulay if and only if $Q:=(p,f,g)\subset S$ is a two generated ideal or all of $S$. Moreover, $p.d._S(R)\leq 1$ and $\nu_S(R)\leq p^2+1$.
    \item If $Q:=(p,f,g)\subseteq S$ has grade three, $R$ admits a birational maximal Cohen-Macaulay module.
\end{enumerate} 
\end{Theorem*}
In section 2, we set up convention and make some preliminary remarks that will be used subsequently. In section 3, we provide some sufficient conditions for $R$ to be Cohen-Macaulay by showing a more general version of \ref{T2}(1), see \ref{T1} and \ref{T1'}. We also identify the conductor of $R$ to the complete intersection ring $A:=S[\omega,\mu]$ in a crucial case. In section 4, we prove parts (2) and (3) of \ref{T2}. We see that $R$ is not ``too far" from being Cohen-Macaulay, in the sense that it can be generated by $p^2+1$ elements over the base ring $S$ and $p.d._S(R)\leq 1$. However, it is not as close as it appears, since if $dim(S)\geq 3$ it could be that $R$ does not even satisfy Serre's condition $S_3$. We then show that $R$ admits a birational maximal Cohen-Macaulay module when $Q:=(p,h_1,h_2)$ has grade three. The condition on $Q$ can be viewed as a further generality condition on the chosen elements. For an $A$-module $M$, let $M^*=Hom_A(M,A)$ be the dual module. The strategy here is to realize $R=I^*$ for a suitable ideal $I\subseteq A$ and identifying an ideal $J\subseteq A$ such that $J^*$ is an $I^*$ module and $J^*$ is $S$-free. 
\section{Preliminaries}
All rings considered are commutative and Noetherian and all modules finitely generated.
\begin{conv}\normalfont\label{notation}
\noindent
\begin{enumerate}
\item Let $R$ be a Noetherian ring and $M$ an $R$-module. For $G\subseteq M$ a subset, the notation $M=\langle G\rangle_R$ means that $M$ is generated as a $R$-module by $G$.
\item Let $(R,\mathfrak{m})$ be a local ring of dimension $d$. A nonzero $R$-module $M$ is a \textbf{maximal Cohen-Macaulay module} (MCM) over $R$ if it is finitely generated and every (some) system of parameters of $R$ is a regular sequence on $M$. If $R$ is an arbitrary Noetherian ring, then an $R$-module $M$ is a \textbf{maximal Cohen-Macaulay module} if for all maximal ideals $\mathfrak{m}\subseteq R$, $M_{\mathfrak{m}}$ is an MCM module over $R_{\mathfrak{m}}$.
\item Let $R$ be a domain and $M$ an $R$-module. Denote by $M^*_R$, the dual module $Hom_R(M,R)$. If $R$ is clear from the context, denote it by $M^*$. Suppose that $M\subseteq K$, where $K$ denotes the field of fractions of $R$. Then via the identification $Hom_R(M,R)\simeq (R:_K M)$, we use $M^*$ to denote $(R:_K M)$ as well.
\item If $R$ is a Noetherian ring of dimension at least one, denote \[NNL_1(R):=\{P\in Spec(R)\: | \: height(P)=1,\: R_P\text{ is not a DVR}\}\]
\item Suppose $S$ is a ring and $p\in \mathbb{Z}$ is such that $p\in S$ is a non-unit. Let $F:S/pS\rightarrow S/pS$ be the Frobenius map. Let $S^p$ denote the subring of $S$ obtained by lifting the image of $F$ to $S$. Let $S^{p^k\wedge p^n}$ for $k,n\geq 1$ denote the multiplicative subset of $S$ of elements expressible in the form $x^{p^k}+yp^n$ for some $x,y\in S$. In particular, $S^{p\wedge p}=S^p$.
\item For a local ring $R$ and an $R$-module $M$, denote $\nu_R(M)$ for the minimal number of generators of $M$ over $R$.
\end{enumerate}
\end{conv}
\begin{rem}\label{R1}\normalfont
Let $S\subseteq C\subseteq D$ be an extension of Noetherian domains such that $S$ is integrally closed, $D$ is module finite over $S$ and $D$ is birational to $C$. Then if $C$ is regular in codimension one, so is $D$. (see \cite{V}[Theorem 2.4] for example)
\end{rem}

\begin{rem}\label{A1}\normalfont
Let $S\subseteq D$ be an extension of Noetherian domains such that going down holds. Let $\overline{D}$ denote the integral closure of $D$ in its field of fractions $K$ and assume $\overline{D}$ is finite over $D$. If $c\cdot u,c\cdot v\in D :_K\overline{D}$ with $c\in D$ and $u.v\in S$ such that there exists no height one prime of $S$ containing both of them, then $NNL_1(D)\subseteq V(c)$.
\end{rem}

\begin{rem}\normalfont\label{P9}
Let $A$ be a Noetherian Gorenstein local domain and $R$ the integral closure of $A$. Assume that $R$ is module finite over $A$. Then for every height one unmixed ideal $I\subseteq A$, $I^*$ is Cohen-Macaulay if and only if $A/I$ is Cohen-Macaulay (see \cite{PS}[Proposition 2.11] for example). In particular, $R$ is Cohen-Macaulay if and only if the conductor of $R$ to $A$ is Cohen-Macaulay. 

\end{rem}
\begin{rem}\label{Gor}\normalfont
Let $S\subseteq R$ be a finite extension of local rings such that $S$ is Gorenstein and $R$ is Cohen-Macaulay. Then for an $R$-module $M$, we have $Hom_R(M,\omega_R)\simeq Hom_S(M,S)$ as $R$-modules (and $S$-modules), where $\omega_R$ is the canonical module of $R$.
\end{rem}
\begin{conv}\normalfont\label{mainconv}
We will assume the following notation for the rest of the paper unless otherwise specified. Let $S$ denote a Noetherian, integrally closed domain of dimension $d$ and $L$ its field of fractions. Assume $Char(L)=0$. Fix $3\leq p\in \mathbb{Z}$ a prime and assume that $p\in S$ is a principal prime such that $S/pS$ is integrally closed. Note here that an unramified regular local ring of mixed characteristic $p$ satisfies the above hypothesis, though not all results here require this specific setting.
\par
Say that a subset $E\subset S$ satisfies $\mathscr{A}_1$, if for all distinct $x,y\in E$, there exists no height one prime $Q\subset S$ such that $(x,y)\subset Q$. An element $c\in S$ is square free if for all height one primes $P\subset S$ containing $c$, $PS_{P}=(c)S_{P}$.
\par
Fix $f,g \in S^p$ square free, non-units, satisfying $\mathscr{A}_1$ or $f,g\in S^p$ units such that they are not $p$-th powers in $S$. Write $f=h_1^p+a\cdot p$ and $g=h_2^p+b\cdot p$ with $h_1,h_2,a,b\in S$.
\par Let $W,U$ be indeterminates over $S$. We have the monic irreducible polynomials $F(W):=W^p-f\in S[W]$ and $G(U):=U^p-g\in S[U]$. Let $\omega$, $\mu$ be roots of $F(W)$ and $G(U)$ respectively and set $K:=L(\omega,\mu)$. Assume that $G(U)$ is irreducible over $L(\omega)$, so that $[K:L]=p^2$. If $S$ is a unique factorization domain and $f,g$ are non units, it can be shown that $[K:L]=p^2$ is automatic. Let $R$ be the integral closure of $S$ in $K$, that is $R$ is the integral closure of $A:=S[\omega,\mu]$. Note that, $A\simeq S[W,U]/(F(W),G(U))$, $S[\omega]\simeq S[W]/(F(W))$ and $S[\mu]\simeq S[U]/(G(U))$. 
\end{conv}
\begin{rem}\normalfont
From \cite{PS}[Proposition 2.10], if exactly one of $f,g$ lies in $S^p$ then $R$ is $S$-free. Moreover \cite{PS}[Example 2.12] shows that $R$ need not be $S$-free when $f,g\notin S^p$. However, to construct an MCM module over $R$ it suffices to consider the case $f,g\in S^p$ when $S$ is a complete unramified regular local ring with perfect residue field, see \cite{DKPS}. This motivates us to understand the case $f,g\in S^p$.
\end{rem}

\begin{Lemma}[\cite{KA}]\label{P6}
Let $p=2k+1$ and $h\in S\setminus pS$. Let $W$ be an indeterminate over $S$. If
\begin{equation}
    C:=(W^p-h^p)-(W-h)^p=\sum_{j=1}^k (-1)^{j+1}{p\choose j}(W\cdot h)^j[W^{p-2j}-h^{p-2j}]
\end{equation}
$C':=(p(W-h))^{-1}\cdot C$ and $\tilde{P}:=(p,W-h)S[W]$, then $C'\notin \tilde{P}$.
\end{Lemma}
\begin{Lemma}\label{c'lemma}
Let $p=2k+1$ and $h\in S\setminus pS$. Let $W$ be an indeterminate over $S$. Suppose $C'$ is as defined in \ref{P6}. Then $C'\equiv h^{p-1}\text{ mod }(p,W-h)S[W]$.
\end{Lemma}
\demo{We have in $S[W]$
 \begin{align*}
 \begin{split}
 C'  &= \sum_{j=1}^k(-1)^{j+1}j^{-1}{p-1\choose j-1}(W\cdot h)^j[W^{p-2j-1}+\dots +h^{p-2j-1}]
 \\
 &\equiv \sum_{j=1}^k(-1)^{j+1}j^{-1}{p-1\choose j-1}h^{2j}\cdot(p-2j)\cdot h^{p-2j-1}\; \text{mod}\: (p,W-h)
 \\
 &\equiv -2h^{p-1}\sum_{j=1}^k(-1)^{j+1}{p-1\choose j-1}\; \text{mod}\: (p,W-h)
 \\
 &\equiv -h^{p-1}(-1)^{k+1}{2k\choose k}\; \text{mod}\: (p,W-h)
 \\
 &\equiv h^{p-1}\; \text{mod}\: (p,W-h)
 \end{split}
 \end{align*}
 \QED
}
\begin{conv}\label{cprime}\normalfont 
Suppose $h_1,h_2\in S\setminus pS$. In this case, let $C_1'$ and $C_2'$ denote respectively the elements in the rings $S[W]$ and $S[U]$ obtained by setting $h=h_1$ and $h=h_2$ in \ref{P6}. Denote by $c_1'$ and $c_2'$ their respective images in the rings $S[\omega]$ and $S[\mu]$ respectively. If $h_1=0$ ($h_2=0$), simply set $c_1'=0$ ($c_2'=0$). Denote by $d_i$ the corresponding element in $S[\omega\mu^i]$ for $1\leq i\leq p-1$.
\end{conv}
\begin{Proposition}[\cite{KA}]\label{P2}
With notation as specified above, $S[\omega]$ is integrally closed if and only if $f\notin S^{p\wedge p^2}$. Further, if $S[\omega]$ is not integrally closed, write $f=h_1^p+a'\cdot p^2$ for some $a'\in S$. Then
\begin{itemize}
\item[{\rm (a)}] $\overline{S[\omega]}=(P^*)_{S[\omega]}=S[\omega,\tau_1]$ where $\tau_1=p^{-1}\cdot(\omega^{p-1}+h_1\omega^{p-2}+\dots+h_1^{p-1})$ and $P:=(p,\omega-h_1)S[\omega]$.
\item[{\rm (b)}] There are exactly two height one primes in $\overline{S[\omega]}$ containing $p$, namely $P:=(p,\omega-h_1,\tau_1)\overline{S[\omega]}$ and $Q:=(p,\omega-h_1,\tau-c_1')\overline{S[\omega]}$. Further, $P_P=(\omega-h_1)_P$ and $Q_Q=(p)_Q$.
\item[{\rm (c)}] The element $\tau_1$ satisfies $l(T):=T^2-c_1'T-a'\cdot(\omega-h_1)^{p-2}\in S[\omega][T]$.
\item[{\rm (d)}] $\overline{S[\omega]}=\langle 1,\omega,\omega^2,\dots,\omega^{p-2},\tau_1\rangle_S$ is $S$-free.
\end{itemize}
\end{Proposition}
\section{Cohen-Macaulay integral closures}
In this section we identify scenarios where $R$ is $S$-free by showing more general versions of $1(a)$ and $1(b)$ in Theorem \ref{T2}. We maintain notation as set up in \ref{mainconv}.
\begin{Proposition} \label{T1} $R$ is $S$-free if at least one of the rings $S[\omega],S[\mu]$ is not integrally closed.
\end{Proposition}
\demo We organize the proof as follows:
\begin{enumerate}
    \item Assume $S[\omega]$ and $S[\mu]$ are both not integrally closed. We then
    \begin{enumerate}
        \item Identify a finite birational overring  $A\hookrightarrow\mathscr{R}A$ such that $\mathscr{R}A$ satisfies $R_1$.
        \item Identify a ``natural" finite birational overring $\mathscr{R}A\hookrightarrow Z$ such that $Z$ is $S$-free, so that $R=Z$ is $S$-free.
    \end{enumerate}
    \item Assume exactly one of the rings $S[\omega]$, $S[\mu]$ is integrally closed. We then take an identical path as indicated in (1) above.
\end{enumerate}
\begin{enumerate}[wide, labelwidth=!, labelindent=0pt]
\item 
\begin{enumerate}[wide, labelwidth=!, labelindent=0pt]
\item From $\ref{P2}$, $f,g\in S^{p\wedge p^2}$. Write $f=h_1^p+a'\cdot p^2$ and $g=h_2^p+b'\cdot p^2$ for some $a',b'\in S$. Note that $h_1,h_2\neq 0$ since $f,g$ are square free. We have from \ref{P2} that $S[\omega,\tau_1], S[\mu,\tau_2]$ are the respective normalizations of $S[\omega]$ and $S[\mu]$ where \[\tau_1=p^{-1}\cdot(\omega^{p-1}+h_1\omega^{p-2}+\dots+h_1^{p-1})\] \[\tau_2=p^{-1}\cdot(\mu^{p-1}+h_2\mu^{p-2}+\dots+h_2^{p-1})\]
\\
Set $E:=A[\tau_1,\tau_2]$. Let $X,Y$ be indeterminates over $A$ and let $\phi:A[X,Y]\rightarrow E$ be the projection map sending $X\shortrightarrow \tau_1$ and $Y\shortrightarrow \tau_2$. From \ref{P2}: \[X^2-c_1'X-a'(\omega-h_1)^{p-2}, Y^2-c_2'Y-b'(\mu-h_2)^{p-2}\in Ker(\phi)\] Height one primes in $E$ containing $p$ correspond to height three primes in $A[X,Y]$ containing $Ker(\phi)$ and $p$. Since $P:=(p,\omega-h_1,\mu-h_2)$ is the unique height one prime in $A$ containing $p$, any such height three prime in $A[X,Y]$ has to contain either $X$ or $X-c_1'$. Likewise it contains either $Y$ or $Y-c_2'$. Therefore if $Q\subseteq A[X,Y]$ is a height three prime containing $p$ and $Ker(\phi)$, it must be that $(p,\omega-h_1,\mu-h_2,X-m,Y-n)\subseteq Q$ for some $m,n\in A$ and hence the containment must be an equality. Moreover, there is at least one height one prime in $E$ containing $p$, since $p$ is not a unit in $S$. Therefore, the only possibilities for height one primes in $E$ containing $p$ are 
\[P_1:=(p,\omega-h_1,\mu-h_2,\tau_1,\tau_2)\] \[P_2:=(p,\omega-h_1,\mu-h_2,\tau_1,\tau_2-c_2')\] \[P_3:=(p,\omega-h_1,\mu-h_2,\tau_1-c_1',\tau_2)\] \[P_4:=(p,\omega-h_1,\mu-h_2,\tau_1-c_1',\tau_2-c_2')\] 
We have $\omega\cdot F'(\omega)=p\cdot f\in (S[\mu,\tau_2,\omega]:_K R)$ and identically $p\cdot g\in (S[\omega,\tau_1,\mu]:_KR)$ (see for example \cite{HS}[Theorem 12.1.1]) and hence $p\cdot f,p\cdot g\in (E:_KR)$. From \ref{A1}, $NNL_1(E)\subseteq V(p)$. But the localizations of $E$ at $P_2,P_3$ and $P_4$ are regular with uniformizing parameters being the images of $\omega-h_1,\mu-h_2$ and $p$ respectively. For example, consider ${P_2}_{P_2}$. Let $Q_1:=(p,\omega-h_1,\tau_1)S[\omega,\tau_1]$ and $Q_2:=(p,\mu-h_2,\tau_2-c_2')S[\mu,\tau_2]$. From $\ref{P2}$, ${Q_1}_{Q_1}=(\omega-h_1)_{Q_1}$ and ${Q_2}_{Q_2}=(p)_{Q_2}$. Thus ${P_2}_{P_2}=(\omega-h_1)_{P_2}$. The $P_3$ and $P_4$ cases are similar. Note that however ${P_1}_{P_1}=(\omega-h_1,\mu-h_2)_{P_1}$.
\par
Set $\eta_1=p^{-1}(\omega-h_1)(\mu-h_2)^{p-2}\in K$. Let $X$ be an indeterminate over $E$. Then $\eta_1$ satisfies $l(X)\in E[X]$ where $l(X):= X^{p-1}-(\tau_1-c_1')(\tau_2-c_2')^{p-2}$ since $p\cdot \tau_1=(\omega-h_1)^{p-1}+p\cdot c_1'$ (similarly for $p\cdot \tau_2$). We claim that $\mathscr{R}A:=E[\eta_1]$ is regular in codimension one. From $\ref{R1}$, $NNL_1(\mathscr{R}A)\subseteq V(P_1\mathscr{R}A)$. Denote by $\overline{l(X)}$ the image of $l(X)$ in $(E/P_1E)[X]$. From \ref{c'lemma} we have $c_1'\equiv h_1^{p-1}$ and $c_2'\equiv h_2^{p-1}$ in the ring $E/P_1E$, and thus \footnote{If $R$ is a ring of characteristic $p$ and $X,Y$ indeterminates over $R$, then for $X^{p-1}-Y^{p-1}\in R[X,Y]$, $X^{p-1}-Y^{p-1}=\prod_{i=1}^{p-1} (X+iY)$.} 
\[\overline{l(X)}= X^{p-1}-(h_1h_2^{p-2})^{p-1}=\prod_{k=1}^{p-1} (X+kh_1h_2^{p-2})\in (E/P_1E)[X]\] Thus, the only possibilities for height one primes in $\mathscr{R}A$ lying over $P_1$ are \[Q_k=(p,\omega-h_1,\mu-h_2,\tau_1,\tau_2,\eta_1+kh_1h_2^{p-2})\mathscr{R}A\] for $1\leq k \leq p-1$ and we have \[{Q_k}_{Q_k}=(\omega-h_1,\mu-h_2,\eta_1+kh_1h_2^{p-2})_{Q_k}\] Since $(\mu-h_2)\eta_1=(\tau_2-c_2')(\omega-h_1)$ and \[\eta_1,\tau_2-c_2', \prod_{j=1,j\neq k}^{p-1}(\eta_1+jh_1h_2^{p-2})\notin Q_k\] we have ${Q_k}_{Q_k} =(\mu-h_2)_{Q_k}=(\omega-h_1)_{Q_k}$. Therefore $\mathscr{R}A$ is regular in codimension one.
\item Set $Z:=\langle T\rangle_S$ where $T:=T_1\cup T_2\cup T_3$ and
\[T_1:=\{(\omega-h_1)^i(\mu-h_2)^j \, | \, 0\leq i,j \leq p-1, i+j<p-1\, \}\]
\[T_2:=\{ p^{-1}(\omega-h_1)^i(\mu-h_2)^j\, | \,  0\leq i,j \leq p-1, i+j\geq p-1,\, i+j\neq 2p-2\, \}\] \[T_3:=\{p^{-2}(\omega-h_1)^{p-1}(\mu-h_2)^{p-1}\}\] For every choice of $0\leq i,j\leq p-1$, there is a unique element in $T$ with ``leading coefficient" $\omega^i\mu^j$. Therefore the order of $T$ is $p^2$. Moreover, since $A$ is $S$-free with a basis given by $D:=\{(\omega-h_1)^i(\mu-h_2)^j \:| \: 0\leq i,j\leq p-1\}$, the elements of $T$ are linearly independent over $S$. From the relations 
\[(\omega-h_1)^p=a'p^2-pc_1'(\omega-h_1)\]
\[(\mu-h_2)^p=b'p^2-pc_2'(\mu-h_2)\]
we get that $Z$ is a ring. Since $Z$ satisfies $S_2$, if we show $\mathscr{R}A\subseteq Z$, then from $\ref{R1}$ $Z=R$. Since $D\subseteq T$, we have $A\subseteq Z$. Since $\eta_1=p^{-1}(\omega-h_1)(\mu-h_2)^{p-2}\in T$, it only remains to be seen that $\tau_1,\tau_2\in Z$. But this is clear from the relation $\tau_1=p^{-1}(\omega-h_1)^{p-1}+c_1'$ (analogously for $\tau_2$). Thus $Z=R$ and $R$ is $S$-free.
\end{enumerate}
\item
\begin{enumerate}[wide, labelwidth=!, labelindent=0pt]
\item Assume without loss of generality $\overline{S[\omega]}=S[\omega,\tau]$ where $\tau = p^{-1}(\omega^{p-1}+\dots+h_1^{p-1})$ and that $S[\mu]$ is integrally closed. Notice that if $P:=(p,\mu-h_2)\subseteq S[\mu]$ is the unique height one prime in $S[\mu]$ containing $p$, then $P_P=(\mu-h_2)_P$ since $(\mu-h_2)(\mu^{p-1}+\dots+h_2^{p-1})=bp$ and $b\notin pS$. From $\ref{P2}$, $f\in S^{p\wedge p^2}$, so write $f=h_1^p+a'p^2$.
\par Set $E:=S[\omega,\mu,\tau]$. From $\ref{P2}$, it follows that there are precisely two height one primes in $E$ containing $p$, namely $P_1:=(p,\omega-h_1,\tau,\mu-h_2)$ and $P_2:=(p,\omega-h_1,\tau-c_1',\mu-h_2)$. From \ref{A1}, $E$ is regular in codimension one outside of $P_1$ and $P_2$. It follows from $\ref{P2}$(b) that ${P_1}_{P_1}=(\omega-h_1,\mu-h_2)_{P_1}$ and ${P_2}_{P_2}=(p,\mu-h_2)_{P_2}=(\mu-h_2)_{P_2}$.
\par Set $\eta_1:=p^{-1}(\omega-h_1)(\mu-h_2)^{p-1}\in K$. $\eta_1\in R$ since it satisfies \begin{equation}\label{latesteqn}
l(X):=X^{p-1}-(\tau-c_1')k_2^{p-2}(\mu-h_2)\in E[X]
\end{equation}
where $k_2=p^{-1}(\mu-h_2)^p\in A\setminus P$. Set $\mathscr{R}A:=E[\eta_1]$. From \ref{R1}, $\mathscr{R}A$ is regular in codimension one if height one primes in $\mathscr{R}A$ lying over $P_1E$ are regular. From the above integral equation for $\eta_1$ over $E$, it is clear that the only such height one prime in $\mathscr{R}A$ is $Q_1:=(p,\omega-h_1,\tau,\mu-h_2,\eta_1)$. Now ${Q_1}_{Q_1}=(\omega-h_1,\mu-h_2,\eta_1)_{Q_1}$. But $\eta_1(\mu-h_2)=k_2(\omega-h_1)$ and $k_2\notin Q_1$. Further, since $\tau-c_1'\notin Q_1$, ${Q_1}_{Q_1}=(\eta_1)_{Q_1}$. Therefore $\mathscr{R}A$ is regular in codimension one.
\item Set $Z:=<T>_S$ where $T=T_1\cup T_2\cup T_3$ and \[T_1:=\{(\omega-h_1)^i(\mu-h_2)^j\;|\; i+j<p, 0\leq i\leq p-2, 0\leq j\leq p-1\}\] \[T_2:=\{p^{-1}(\omega-h_1)^i(\mu-h_2)^j\;|\;i+j\geq p, 1\leq i\leq p-1,1\leq j\leq p-1\} \] \[T_3:=\{\tau-c_1'=p^{-1}(\omega-h_1)^{p-1}\}\] For every $0\leq i,j\leq p-1$, there is a unique element in $T$ with ``leading coefficient" $\omega^i\mu^j$. Therefore the order of $T$ is $p^2$. Moreover since $A$ is $S$-free with a basis given by $D:=\{(\omega-h_1)^i(\mu-h_2)^j \:| \: 0\leq i,j\leq p-1\}$, the elements of $T$ are linearly independent over $S$. From the relations
\[(\omega-h_1)^p=a'p^2-pc_1'(\omega-h_1)\]
\[(\mu-h_2)^p=bp-pc_2'(\mu-h_2)\]
we see that $Z$ is a ring. Since $Z$ satisfies $S_2$, if we show that $\mathscr{R}A\subseteq Z$, then from \ref{R1} $Z=R$. It now suffices to note that $D\subseteq Z$, $\eta_1\in T$ and $ \tau:=p^{-1}(\omega-h_1)^{p-1}+c_1'\in Z$, so that $\mathscr{R}A\subseteq Z$. Thus $R=Z$ is $S$-free. 
\end{enumerate}
\end{enumerate}\QED
\begin{Lemma}\label{lem1}
With established notation, assume that $S[\omega]$ and $S[\mu]$ are integrally closed. The following hold
\begin{enumerate}
    \item $(P^{(p-1)})^*_A=\langle T\rangle_S$, where $P^{(p-1)}$ denotes the $(p-1)$-th symbolic power of the unique height one prime $P\subseteq A$ containing $p$ and $T:=T_1\cup T_2$, with \[T_1:=\{(\omega-h_1)^i(\mu-h_2)^j|\;0\leq i,j\leq p-1\;,\;i+j<p\}\]
\[T_2:=\{p^{-1}(\omega-h_1)^i(\mu-h_2)^{j}|\;0\leq i,j\leq p-1\;,\;i+j\geq p\}\]
\item The ring $A/P^{(p-1)}$ is Cohen-Macaulay.
\end{enumerate}
\end{Lemma}
\demo For every $0\leq i,j\leq p-1$, there is a unique element in $T$ with ``leading coefficient" $\omega^i\mu^j$. Therefore the order of $T$ is $p^2$. Moreover since $A$ is $S$-free with a basis given by $D:=\{(\omega-h_1)^i(\mu-h_2)^j \:| \: 0\leq i,j\leq p-1\}$, the elements of $T$ are linearly independent over $S$. From the relations
\[(\omega-h_1)^p=ap-pc_1'(\omega-h_1)\]
\[(\mu-h_2)^p=bp-pc_2'(\mu-h_2)\]
we see that $\langle T \rangle_S$ is a ring (in particular it is an $A$-module). Moreover, since it is Cohen-Macaulay, (2) immediately follows from (1) by using \ref{P9}. Therefore, only (1) remains to be shown. Since $(P^{(p-1)})^*$ and $\langle T\rangle_S$ are birational, $S_2$ $A$-modules, it suffices to show their equality in codimension one. If $Q\subseteq A$, $Q\neq P$ is a height one prime, the equality is clear. So localize $A$ at $P$ and assume $(A,P)$ local for the rest of the proof. Then $\langle T\rangle_S=A[\eta]$ where $\eta:=(\omega-h_1)^{-1}(\mu-h_2)$.
Note that $A= B'/(G(U))$, where $B':=S[\omega][U]$. Set $\tilde{I}:=(\omega-h_1,U-h_2)^{p-1}\subseteq B'$. We have $G(U)\in \tilde{I}$: \begin{align*}
    U^p-g &=(U-h_2)^p + p(C_2'(U-h_2)-b)
    \\
    &=-k_1^{-1}(\omega-h_1)(C_2'(U-h_2)-b)\cdot \Delta_1 +0\cdot\Delta_2 +\dots+ 0\cdot\Delta_{p-1} - (U-h_2)\cdot\Delta_p 
  \end{align*} where $\Delta_i=(-1)^ie_i$ with $e_i=(\omega-h_1)^{p-i}(U-h_2)^{i-1}$ and $k_1=p^{-1}(\omega-h_1)^p$. That is $\tilde{I}$ is the lift to $B'$ of $P^{p-1}$. Further, $\tilde{I}$ is grade two perfect since it arises as the ideal of maximal minors of the $p\times(p-1)$ matrix $M$ :
  \[
M = \begin{bmatrix} 
    U-h_2 & 0 & \dots & 0 & 0 \\
    \omega-h_1 & U-h_2 & 0 &\dots & 0\\
    0  & \omega-h_1  &U-h_2  &\dots & 0  \\
    \vdots &0 &\ddots & \ddots& \vdots\\
    0 & \dots & 0 & \omega-h_1 & U-h_2 \\
    0 & \dots  & 0 & 0 & \omega-h_1 
    \end{bmatrix}
    \]
 Let $M'$ be the $p\times p$ matrix obtained by adjoining $M$ with the column of coefficients of $G(U)$:
 \[
M' = \begin{bmatrix} 
    U-h_2 & 0 & \dots & 0 & 0& -k_1^{-1}(\omega-h_1)(C_2'(U-h_2)-b)\\
    \omega-h_1 & U-h_2 & 0 &\dots & 0& 0\\
    0  & \omega-h_1  &U-h_2  &\dots & 0 & 0 \\
    \vdots &0 &\ddots & \ddots& \vdots & \vdots\\
    0 & \dots & 0 & \omega-h_1 & U-h_2&0 \\
    0 & \dots  & 0 & 0 & \omega-h_1 & -(U-h_2)
    \end{bmatrix}
\]
By \cite{KU}[Lemma 2.5] (or \cite{KA}[Prop 2.1]), $(P^{p-1})^*$ is generated as an $A$-module by the set $\{\delta_i^{-1}M'_{i,i}\:|\: 1\leq i\leq p\}$, where $\delta_i$ denotes the image of $\Delta_i$ in $A$ and $M'_{i,i}$ the image in $A$ of the $(i,i)$-th cofactor of $M'$. This is exactly the set $\{\eta^{p-1}, \eta^{p-2},\dots,\eta,1\}$. Since $\eta$ satisfies the integral equation $X^p-k_1^{-1}k_2\in A[X]$ (with $k_2=p^{-1}(\mu-h_2)^p$), this implies $(P^{p-1})^*=A[\eta]=\langle T\rangle_S$. Thus the proof is complete. 
\QED
\begin{Proposition}\label{T1'}
With established notation, $R$ is $S$-free if $S[\omega]$ and $S[\mu]$ are integrally closed and $fg^i\notin S^{p\wedge p^2}$ for $1\leq i\leq p-1$.
Further in this case, $P^{(p-1)}$ is the conductor of $R$ to $A$ where $P$ is the unique height one prime in $A$ containing $p$ and $P^{(p-1)}$ denotes the $(p-1)$-th symbolic power of $P$.
\end{Proposition}
\demo Since $S[\omega]$ and $S[\mu]$ are integrally closed, we have from $\ref{P2}$ that $f,g\notin S^{p\wedge p^2}$. Write $f=h_1^p+ap$ and $g=h_2^p+bp$ with $a,b\notin pS$. We first note the following: The condition $fg^i\notin S^{p\wedge p^2}$ for all $1\leq i\leq p-1$ is equivalent to the condition $\prod_{i=1}^{p-1}(ah_2^p+ibh_1^p)\notin pS$. This follows since for $1\leq i\leq p-1$
\begin{equation*}
\begin{aligned}
fg^i&=(h_1^p+ap)(h_2^p+bp)^i
\\
&=(h_1h_2^i)^p+h_2^{p(i-1)}(ah_2^p+ibh_1^p)\cdot p +q\cdot p^2
\end{aligned}
\end{equation*} 
for some $q\in S$. We organize the proof as follows:
\begin{enumerate}
    \item We first construct the normalization of $A$ locally at $NNL_1(A)$. Since $NNL_1(A)=\{P\}$, we only need to construct $R_P$. 
    \item Using (1), identify a finite birational overring $A\hookrightarrow \mathscr{R}A$ such that $\mathscr{R}A$ satisfies $R_1$. Then choose a suitable finite birational overring $\mathscr{R}A\hookrightarrow Z$ such that $Z$ is $S$-free. This would show $R=Z$ is $S$-free. 
\end{enumerate}
\begin{enumerate}[wide, labelwidth=!, labelindent=0pt]
\item From $\ref{A1}$, $A$ is regular in codimension one outside of $V(p)$. Moreover $P:=(p,\omega-h_1,\mu-h_2)$ is the only height one prime in $A$ containing $p$. Localize at $P$ and assume $(A,P)$ local for part (1). Now $(\omega-h_1)^p=pk_1$ and $(\mu-h_2)^p=pk_2$ where $k_1=a-c_1'(\omega-h_1)$ and $k_2=b-c_2'(\mu-h_2)$. Since $k_1,k_2\notin P$, we have $P=(\omega-h_1,\mu-h_2)$.
 The element $\eta = (\omega-h_1)^{-1}(\mu-h_2)\in K$ satisfies the integral equation $l(X):=X^p-k_1^{-1}k_2\in A[X]$. If $\eta\in A$, then $P_P=(\omega-h_1)_P$ so $A$ is integrally closed. But from \cite{PS}[Proposition 2.7] this is impossible. Therefore $A[\eta]$ is a proper birational extension of $A$. We claim that $E:=A[\eta]$ is regular when $\prod_{i=1}^{p-1}(ah_2^p+ibh_1^p)\notin pS$. We will observe that $E$ is local with maximal ideal $Q\subseteq E$ either of the form $Q=PE$ or $Q=(P,\eta-r)E$ for some suitable $r\in S\setminus pS$. To see this, let $\phi:\tilde{E}:=S[W,U]_{\tilde{P}}[X]\longrightarrow E$ be the natural projection map sending $W\mapsto \omega, U\mapsto \mu, X\mapsto \eta$, where $W,U,X$ are indeterminates over $S$ and $\tilde{P}:=(p,W-h_1,U-h_2)$. Let $Q\subseteq E$ be any maximal ideal and let $\tilde{Q}$ be the preimage of $Q$ under $\phi$.  Now 
 \begin{equation}\label{integraletacharp}
 l(X)\equiv X^p-ba^{-1} \in (A/P)[X]\simeq (S/pS)[X]
 \end{equation}
 Since $S/pS$ is a field, if $l(X)$ is an irreducible polynomial over $(S/pS)[X]$ then \[(p,W-h_1,U-h_2,l(X))\subseteq \tilde{Q}\] is a height four prime containing $p$. Therefore the above inclusion must be an equality.
 If on the other hand $l(X)$ is reducible over $S/pS$ then $l(X)\equiv (X-r)^p \in (S/pS)[X]$ for some $r\in S\setminus pS$. So in this case \[(p,W-h_1,U-h_2,X-r)\subset \tilde{Q}\] is a height four prime. Again, the above inclusion must then be an equality. Therefore in either case $E$ is local and the maximal ideal is either of the form $PE=(\omega-h_1,\mu-h_2)E$ or $Q:=(\omega-h_1, \mu-h_2, \eta-r)E$. In the first case, since $\eta\cdot (\omega-h_1)=\mu-h_2$, $E$ is a DVR. In the second case, we have $Q_{Q}=(\omega-h_1,\eta-r)_Q$. We now show that $Q_Q=(\eta-r)_Q$. We have for some $m\in E$:
 \begin{align*}
 (\eta-r)^p  &= \eta^p -ba^{-1} +pm
 \\
 &= k_1^{-1}k_2-ba^{-1} +pm
 \\
 & = k_1^{-1}[-c_2'(\mu-h_2)+ba^{-1}c_1'(\omega-h_1)+pmk_1]
 \\
 &= k_1^{-1}(\omega-h_1)[-c_2'\eta + ba^{-1}c_1'+(\omega-h_1)^{p-1}m]    
 \end{align*}
 So $Q$ is principal if $\alpha:=c_2'\eta-ba^{-1}c_1'$ is invertible in $E$. To show $\alpha\in E$ is a unit, it suffices to show $(ac_2'\eta-bc_1')^p\in E$ is invertible. From \ref{c'lemma}, $c_1'\equiv h_1^{p-1}\:\text{mod}\: Q$ and $c_2'\equiv h_2^{p-1}\:\text{mod}\: Q$. We then have
  \begin{align*}
 (ac_2'\eta-bc_1')^p  &\equiv a^p(c_2')^p\eta^p-b^p(c_1')^p\:\text{mod}\:Q
 \\
 &\equiv b(a^{p-1}(c_2')^p-b^{p-1}(c_1')^p)\:\text{mod}\:Q
 \\
 & \equiv b[a^{p-1}(h_2^{p-1})^p-b^{p-1}(h_1^{p-1})^p]\:\text{mod}\:Q
 \\
 &\equiv b\prod_{i=1}^{p-1}(ah_2^p+ibh_1^p)\:\text{mod}\:Q  \end{align*}
 Thus $\alpha\in E$ is a unit. Hence $E=A[\eta]$ is regular, that is $E=R$.

\item Set $\mathscr{R}A:=A[k_1\eta]$ for $\eta=(\omega-h_1)^{-1}(\mu-h_2)$ and $k_1=p^{-1}(\omega-h_1)^p$. Note that $k_1\eta=p^{-1}(\omega-h_1)^{p-1}(\mu-h_2)$ and that it satisfies the integral equation $X^p-k_1^{p-1}k_2\in A[X]$ for $k_2:=p^{-1}(\mu-h_2)^p$. Since $k_1\notin P$, by \ref{R1} and part (1) of the proof, $\mathscr{R}A$ is regular in codimension one.
\par Set $Z:=\langle T\rangle_S$ where $T$ is as in the statement of \ref{lem1}. We see that $Z$ is a ring from the relations
\[(\omega-h_1)^p=ap-pc_1'(\omega-h_1)\]
\[(\mu-h_2)^p=bp-pc_2'(\mu-h_2)\]
Moreover it is a free $S$-module of rank $p^2$. Clearly $\mathscr{R}A\subseteq Z$, so $Z$ inherits $R_1$ from $\mathscr{R}A$. Thus $Z=R$ and $R$ is $S$-free. 
\par Finally, from \ref{lem1}(1) $P^{(p-1)}$ is contained in the conductor $J$ of $R$ to $A$. Since $A_P$ is a one dimensional Gorenstein local ring, $J_P=(P^{p-1})_P$ and thus $J\subseteq P^{(p-1)}$. Thus $P^{(p-1)}$ is the conductor of $R$ to $A$.

\end{enumerate}
\QED
\begin{rem} \normalfont
The condition $fg^i\notin S^{p\wedge p^2}$ for $1\leq i\leq p-1$ in $\ref{T1'}$ is saying that some suitable subrings of $A$ are integrally closed. As noted in the proof of $\ref{T1'}$, the condition is equivalent to $\prod_{i=1}^{p-1}(ah_2^p+ibh_1^p)\notin pS$.  Let $1\leq k,i\leq p-1$ and $1\leq i(k)\leq p-1$ be such that $i(k)-ik\in p\mathbb{Z}$. The condition $\prod_{i=1}^{p-1}(ah_2^p+ibh_1^p)\notin pS$ is saying that for all $1\leq i\leq p-1$, $A_i:=S[\omega\mu^{i(1)},\dots,\omega^j\mu^{i(j)},\dots,\omega^{p-1}\mu^{i(p-1)}]$ is integrally closed. Indeed \[(\omega\mu^i)^p=fg^i=(h_1h_2^i)^p+(ah_2^{ip}+ibh_1^ph_2^{ip-p})p+p^2q\] for some $q\in S$. If $i=1$, we have that $fg$ is squarefree in $S$ and by \ref{P2} $S[\omega\mu]$ is integrally closed. If $i\neq 1$, the given condition is equivalent to saying that \[NNL_1(S[\omega\mu^i])\cap V(p)=\emptyset\] Moreover $NNL_1(S[\omega\mu^i])\subseteq V(g)$. Choose $k$ such that $i(k)=1$, so that $S[\omega\mu^i,\omega^k\mu]$ is a finite birational extension of both $S[\omega\mu^i]$ and $S[\omega^k\mu]$. Now $NNL_1(S[\omega^k\mu])\subseteq V(f)$. So by \ref{R1}, $S[\omega\mu^i,\omega^k\mu]$ is regular in codimension one since $f,g\in S$ satisfy $\mathscr{A}_1$. Since $A_i$ is a finite birational extension of $S[\omega\mu^i,\omega^k\mu]$, it is regular in codimension one for the same reason. The remark follows since it is easily checked that $A_i$ is $S$-free.
\end{rem}
\begin{rem}\normalfont
The powers of the prime $P^{p-1}\subseteq A$ in \ref{T1'} are not $P$-primary in general. For example if $p=3$, observe that $3a,3b\in P^2$. However, it holds that $P^{(p-1)}=(p)+P^{p-1}$. 
\end{rem}

\section{Existence of birational MCM modules}

In this section we look at cases where $R$ is not $S$-free. That is in the primary case of interest, when $S$ is an unramified regular local ring of mixed characteristic $p$, we look at non-Cohen-Macaulay integral closures $R$. More specifically, we will show that $p.d._S(R)=1$ under some natural conditions and show that in this case $R$ admits a birational MCM module.
\par From \ref{T1} and \ref{T1'}, if we are looking for a non $S$-free $R$, we must have that $S[\omega]$ and $S[\mu]$ are integrally closed such that there exists an $1\leq i\leq p-1$ satisfying $fg^i\in S^{p\wedge p^2}$. The reader can easily see that if it exists, such an ``$i$" is unique. We start by identifying an ideal $I\subseteq A$, such that $I^*=R$ under this circumstance.
\begin{conv}\normalfont
We maintain notation as set up in \ref{mainconv} and make the additional assumption that $f,g \notin S^{p\wedge p^2}$. Write $f=h_1^p+ap$, $g=h_2^p+bp$ with $a,b\notin pS$. Assume further that $h_1,h_2\neq 0$. Note here that if $h_1=0$ (or $h_2=0$), we have by $\ref{T1'}$ that $R$ is $S$-free. Let $P:=(p,\omega-h_1,\mu-h_2)\subseteq A$ denote the unique height one prime in $A$ containing $p$.
\end{conv}

\begin{Lemma}\label{lem2}
For $H:=(p,\omega\mu^i-h_1h_2^i)\subseteq A$, $H^*_P=\langle 1,\tau_i\rangle_{A_P}$ where 
\[\tau_i=p^{-1}[(\omega\mu^i)^{p-1}+h_1h_2^i(\omega\mu^i)^{p-2}+\dots+(h_1h_2^i)^{p-1}]\in K\]
\end{Lemma}
\demo Localize $A$ at $P$ and assume $(A,P)$ local. Consider the ideal \[\tilde{H}:=(p,W\mu^i-h_1h_2^i)\subseteq S[\mu][W]\] We have $F(W)\in \tilde{H}$ :
\begin{equation}\label{G(U)}
F(W)-h_2^{-ip}[(W\mu^i)^{p-1}+h_1h_2^i(W\mu^i)^{p-2}+\dots +(h_1h_2^i)^{p-1}]\cdot (W\mu^i-h_1h_2^i)\in pS[\mu][W]
    \end{equation}
    
    Clearly $\tilde{H}$ is a grade two perfect ideal in $S[\mu][W]$ and is the ideal of maximal minors of the matrix $E$:
   \[   E = \begin{bmatrix} 
    W\mu^i-h_1h_2^i \\
    p\\
 \end{bmatrix}
    \]
    Adjoining the column of coefficients from ($\ref{G(U)}$) appropriately, we have for some $\alpha\in S[\mu[[W]$ the matrix $E'$:
    \[   E' = \begin{bmatrix} 
   W\mu^i-h_1h_2^i  & \alpha\\
    p& h_2^{-ip}[(W\mu^i)^{p-1}+h_1h_2^i(W\mu^i)^{p-2}+\dots+ (h_1h_2^i)^{p-1}]\\

    \end{bmatrix}
    \]
From \cite{KA}[Proposition 2.1], $H^*=\langle E'_{11}/\delta_1,E'_{22}/\delta_2\rangle_A$ where $E'_{ii}$ and $\delta_i$ denote the image in $A$ of the $(i,i)$-th cofactor of $E'$ and the $i$-th (signed) minor of $E$ respectively. Thus $H^*=\langle 1,\tau_i\rangle_A$ 
\QED

\begin{Lemma}\label{P10}
With established notation, let $fg^i\in S^{p\wedge p^2}$. Then for \[I:=pA+P^{p-2}\cdot(\omega\mu^i-h_1h_2^i)A\] we have $I^*_A=R$. 
\end{Lemma}
\demo 
Since $I^*_A$ and $R$ are birational $S_2$ $A$-modules, it suffices to show the desired equality in codimension one. If $Q\neq P$ is a height one prime in $A$, $I^*_Q=R_Q=A_Q$. Therefore localize $A$ at $P$ and assume $(A,P)$ and $(S,pS)$ are one dimensional local rings for the remainder of the proof.
\par 
We have $A=S[\mu,\omega\mu^i]$. Note that $g,fg^i\in S$ are units and therefore trivially are square free and satisfy $\mathscr{A}_1$. Moreover, $S[\mu]$ is integrally closed and $S[\omega\mu^i]$ is not. Thus we are in the setting of \ref{T1}(2). From the proof of \ref{T1}(2)(a), we get that $R=A[\tau,\eta]$ where $\tau=p^{-1}[(\omega\mu^i)^{p-1}+h_1h_2^i(\omega\mu^i)^{p-2}+\dots+(h_1h_2^i)^{p-1}]$ and $\eta=p^{-1}(\omega\mu^i-h_1h_2^i)(\mu-h_2)^{p-1}$. 
Since $\tau\eta\in A$, we see from \ref{P2}(c) and equation (\ref{latesteqn}) that $R=\langle 1,\eta,\dots,\eta^{p-2},\tau\rangle_A$.
\par Since $p\in P^{p-1}$, a straightforward calculation gives \[P^{p-1}\cap(p,\omega\mu^i-h_1h_2^i)=(p)+P^{p-1}\cap(\omega\mu^i-h_1h_2^i)=(p)+P^{p-2}\cdot (\omega\mu^i-h_1h_2^i)=I\]  From \ref{lem1}, $(P^{p-1})^*=A[\eta]\subseteq R$. Since $A$ is Gorenstein, $A:_K R\subseteq P^{p-1}$. Let $H$ be as in \ref{lem2}. Combining \ref{lem2} and \ref{P2}(c), we get that $H^*=A[\tau]\subseteq R$. Again since $H$ is reflexive, $A:_K R\subseteq H$. Therefore, $A:_K R\subseteq P^{p-1}\cap H=I$.
\par To show $IR\subseteq A$, note that $I\eta^i\subseteq A$ for $1\leq i\leq p-2$, since $I\subseteq P^{p-1}=A[\eta]^*$. Similarly, $I\tau\subseteq A$ since $I\subseteq H=A[\tau]^*$. Thus we have shown $I=A:_K R$ and the proof is complete.
\QED

We now set out to show $R$ need not be Cohen-Macaulay.
\begin{Lemma}\label{lemRgen}
Let $S[\omega],S[\mu]$ be integrally closed such that $fg\in S^{p\wedge p^2}$. Then 
\begin{enumerate}
    \item $R\subseteq \langle\{1\}\cup p^{-1}\cdot (\omega-h_1,\mu-h_2)^{p-1}\rangle_A$.
    \item  Consider $y=p^{-1}(\sum_{i=1}^p a_i(\mu-h_2)^{p-i}(\omega-h_1)^{i-1})\in K$
with the $a_i\in A$. Then $y\in R$ if and only if for all $2\leq i\leq p$, $a_{i-1}h_2+a_ih_1\in P$.
\end{enumerate} 
\end{Lemma}
\demo 
From \ref{P10}, $p\cdot R\subseteq A$, so consider an arbitrary element $y:=p^{-1}\cdot x\in R$ with $x\in A$. From \ref{P10}, $x\cdot (\omega-h_1)^{p-2}(\omega\mu-h_1h_2)\in pA$. Lifting to $B:=S[W,U]$ and denoting lifts by $\sim$
\begin{equation}
    \tilde{x}(W-h_1)^{p-2}(WU-h_1h_2) \in (p,F(W),G(U))
\end{equation}
Write
\begin{equation}\label{omegamui}
\omega\mu-h_1h_2=(\omega-h_1)(\mu-h_2)+h_2(\omega-h_1)+h_1(\mu-h_2)
\end{equation}
Lifting the identity in (\ref{omegamui}) to $B$ we see that 
$\tilde{x}\in (p,W-h_1,(U-h_2)^{p-1})$. By symmetry $\tilde{x}\in (p,(W-h_1)^{p-1},U-h_2)$ and hence \[\tilde{x}\in (p,(U-h_2)^{p-1},(W-h_1)^{p-1},(W-h_1)(U-h_2))\] This is because for a regular sequence $(q,y,z)\subseteq B$ \begin{equation}\label{regsequence}
(q,y,z^n)\cap (q,y^n,z)=(q,y^n,z^n,yz)
\end{equation}
Since $1\in R$, towards describing $A$-module generators for $R$ we may assume that $y=p^{-1}x$ with 
\begin{equation}\label{Amod1}
x=a_1\cdot (\mu-h_2)^{p-1}+a_2\cdot (\omega-h_1)^{p-1}+a_3\cdot (\omega-h_1)(\mu-h_2)
\end{equation}
for some $a_1,a_2,a_3\in A$. Suppose we can write \begin{equation}\label{Amod2}
y=p^{-1}[a_1(\mu-h_2)^{p-1}+a_2(\omega-h_1)(\mu-h_2)^{p-2}+\dots+a_p(\omega-h_1)^{p-1}+b\cdot(\omega-h_1)^i(\mu-h_2)^i]
\end{equation}
with $1\leq i< (p-1)/2$ and $a_i,b\in A$. By (\ref{Amod1}), we can do this for $i=1$, where $a_j=0$ for $2\leq j\leq p-1$. Now using (\ref{omegamui}) we get that $y\cdot (\omega-h_1)^{p-1-i}(\mu-h_2)^{i-1}(\omega\mu-h_1h_2)\in pA$ if and only if 
\begin{align*}
 a_{i+1}h_1(\omega-h_1)^{p-1}(\mu-h_2)^{p-1}+a_ih_2(\omega-h_1)^{p-1}(\mu-h_2)^{p-1}+bh_1(\omega-h_1)^{p-1}(\mu-h_2)^{2i}\in pA
\end{align*}
Pulling back to $B$
\begin{equation}
\begin{aligned}
    &\quad \quad \tilde{b}h_1(U-h_2)^{2i}+h_2\tilde{a_i}(U-h_2)^{p-1}+h_1\tilde{a}_{i+1}(U-h_2)^{p-1} \in (p,(U-h_2)^p,W-h_1)
    \end{aligned}
    \end{equation}
    and thus $\tilde{b}\in (p,W-h_1,(U-h_2)^{p-2i-1})$. By symmetry, $\tilde{b}\in (p,U-h_2,(W-h_1)^{p-2i-1})$ and hence by (\ref{regsequence}) 
\[\tilde{b}\in (p,(U-h_2)^{p-2i-1},(W-h_1)^{p-2i-1},(W-h_1)(U-h_2))\] Therefore \[py\in (p,(\omega-h_1)^{i+1}(\mu-h_2)^{i+1})+(\omega-h_1,\mu-h_2)^{p-1}\] Starting from (\ref{Amod1}) and iterating the argument from (\ref{Amod2}) to this point sufficiently many times, we see that \[R\subseteq \langle\{1\}\cup p^{-1}\cdot (\omega-h_1,\mu-h_2)^{p-1}\rangle_A\] Consider
\begin{equation}\label{Amod3}
    y=p^{-1}(\sum_{i=1}^p a_i(\mu-h_2)^{p-i}(\omega-h_1)^{i-1})\in K
\end{equation}
with the $a_i\in A$. From \ref{P10}, $y\in R$ if and only if for all $2\leq i\leq p$
\begin{equation}\label{Amod4}
    y\cdot (\omega-h_1)^{p-i}(\mu-h_2)^{i-2}(\omega\mu-h_1h_2)\in A
    \end{equation}
From (\ref{omegamui}) the above statements are equivalent to
\begin{equation}\label{equivalence}
    (a_{i-1}h_2+a_ih_1)(\mu-h_2)^{p-1}(\omega-h_1)^{p-1} \in pA
\end{equation}
for each $2\leq i\leq p$. Lifting to $B$, we see that (\ref{equivalence}) is equivalent to
\begin{equation}\label{equiv1}
a_{i-1}h_2+a_ih_1\in P
\end{equation}
Thus the proof is complete.
\QED

\begin{Lemma}\label{lem2Rgen}
Assume $S$ is an unramified regular local ring of mixed characteristic $p\geq 3$. Let $S[\omega],S[\mu]$ be integrally closed such that $fg\in S^{p\wedge p^2}$. Then $\nu_S(R)\leq p^2+1$. More explicitly, set $\eta_i:=p^{-1}(\omega-h_1)^i(\mu-h_2)^{p-i}$ for $1\leq i\leq p-1$. We have $R=\langle A[\eta_1,\dots,\eta_{p-1}]\cup\{\epsilon\}\rangle_S$ for:
\[\epsilon:=p^{-1} \sum_{i=1}^p (-1)^i c^{p-i}e^{i-1}(\mu-h_2)^{p-i}(\omega-h_1)^{i-1}\]
where $h_1\equiv zc\:\text{mod}\:pS$, $h_2\equiv ze\:\text{mod}\:pS$ for some $z\in S\setminus pS$ and $c,e\in S$ relatively prime.
\end{Lemma}
\demo Suppose $\langle T\rangle_S$ is as in \ref{lem1}. Note that it is just the ring $A[\eta_1,\dots,\eta_{p-1}]$. Since $\langle T\rangle_S$ is $S$-free of rank $p^2$, the assertion $\nu_S(R)\leq p^2+1$ follows from the second assertion. 
\par
Since $S/pS$ is regular local (a UFD), $h_1\equiv (zc)\:\text{mod}\:pS$, $h_2\equiv (ze)\:\text{mod}\:pS$ for some $z\in S\setminus pS$ and $c,e\in S$ relatively prime. First suppose $c$ or $e$ is a unit in $S$. Then it follows from \ref{lemRgen}(1) and (2) that $R=\langle A[\eta_1,\dots,\eta_{p-1}]\cup {\epsilon}\rangle_A$. Notice that if $i+j>0$, then $(\omega-h_1)^i(\mu-h_2)^j\in (A[\eta_1,\dots,\eta_{p-1}]:_K \epsilon)$. Thus $R=\langle A[\eta_1,\dots,\eta_{p-1}]\cup \{\epsilon\}\rangle_S$ in this case.
\par Next assume that neither $c$ or $e$ is a unit, so that $(p,c,e)\subseteq S$ forms a regular sequence. Now $\epsilon\in R$ from \ref{lemRgen}(2). From \ref{lemRgen}(1), it suffices to look at elements of the form \[y=p^{-1}(\sum_{i=1}^p a_i(\mu-h_2)^{p-i}(\omega-h_1)^{i-1})\in R\] In view of \ref{lemRgen}(2), the condition $a_1h_2+a_{2}h_1\in P$ upon lifting to $B:=S[W,U]$ (denoting lifts by $\sim$) tells us that $\tilde{a_1},\tilde{a_2}$ arise from the first syzygy of the grade five complete intersection $B$-ideal, $\tilde{Q}:=(p,c,e,W-h_1,U-h_2)$. In particular \[\tilde{a_2}\in (p,c,W-h_1,U-h_2)\cap (p,e,W-h_1,U-h_2)=(p,ce,W-h_1,U-h_2)\]  since $a_2h_2+a_3h_1\in P$ as well. Since $A[\eta_1,\dots,\eta_{p-1}]\subseteq R$, towards describing $A$-module generators for $R$ we may assume $\tilde{a_2}=\alpha ce$ for some $\alpha\in B$ and consequently that $\tilde{a_1}=-\alpha c^2$ and $\tilde{a_3}=-\alpha e^2$. 
\par
Now let $3\leq i<p$ be such that for all $1\leq k\leq i$, $a_k=(-1)^k\alpha c^{i-k}e^{k-1}$ for some $\alpha\in A$. Lifting $a_{i}h_2+a_{i+1}h_1\in P$ to $B$, we have $(-1)^i\alpha e^i+a_{i+1}c\in \tilde{P}$. Since $A[\eta_1,\dots,\eta_{p-1}]\subseteq R$, we may assume that $\alpha = \alpha'c$ for some $\alpha'\in B$ and hence that $a_{i+1}=(-1)^{i+1}\alpha'\cdot e^i$. Iterating this argument, we get that $R=\langle A[\eta_1,\dots,\eta_{p-1}]\cup\{\epsilon\}\rangle_A$. Finally, if $i+j>0$ then $(\omega-h_1)^i(\mu-h_2)^j\in (A[\eta_1,\dots,\eta_{p-1}]:_K \epsilon)$ and the conclusion follows. 
\QED
\begin{Proposition}\label{RnotCM}
Let $(S,\mathfrak{m})$ be an unramified regular local ring of mixed characteristic $p\geq 3$ such that $S[\omega]$ and $S[\mu]$ are integrally closed and $fg\in S^{p\wedge p^2}$. Then $R$ is Cohen-Macaulay if and only if $Q:=(p,h_1,h_2)\subset S$ is a two generated ideal or all of $S$. Moreover, $p.d_S(R)\leq 1$.
\end{Proposition}
\demo
Since $S/pS$ is a UFD, $h_1\equiv zc\:\text{mod}\:pS$, $h_2\equiv ze\:\text{mod}\:pS$ for some $z\in S\setminus pS$ and $c,e\in S$ relatively prime. Then $Q$ is a two generated ideal or all of $R$ if and only if $c$ or $e$ is a unit. From \ref{lem2Rgen}, $R=\langle A[\eta_1,\dots,\eta_{p-1}]\cup \{\epsilon\}\rangle_S$. Suppose that $c$ is a unit. Then
\[(\mu-h_2)^{p-1}\in \langle \epsilon,(\omega-h_1)(\mu-h_2)^{p-2},\dots,(\omega-h_1)^{p-1}\rangle_S.\]  Thus $R$ is $S$-free of rank $p^2$ and hence is Cohen-Macaulay.
\par Now assume neither $c$ nor $e$ is a unit, that is $Q$ is either grade three perfect or grade two and not perfect. We know from \ref{lem2Rgen} that $R=\langle A[\eta_1,\dots,\eta_{p-1}]\cup\{\epsilon\}\rangle_S$. With $T$ as in \ref{lem1}, define $\Gamma:T\rightarrow \mathbb{Z}$, $\Gamma':T\rightarrow \mathbb{Z}$ by \[\Gamma(p^{-k}(\omega-h_1)^i(\mu-h_2)^j)=i+j\]
\[\Gamma'(p^{-k}(\omega-h_1)^i(\mu-h_2)^j)=i\]
Define a total ordering on $T$ as follows: for $x,y\in T$, if $\Gamma(x)\geq \Gamma(y)$ then $x\geq y$ and if $\Gamma(x)=\Gamma(y)$, then $x\geq y$ if $\Gamma'(x)\geq \Gamma'(y)$. Let $\alpha:S^{p^2+1}\rightarrow R$ be the $S$-projection map defined by the generating set $T\cup \{\epsilon\}$ such that the basis element $e_{p^2+1}$ maps to $\epsilon$ and the image of the basis elements $e_i$, $i\neq p^2+1$ is defined by the ordered set $T$. Consider $U=[u_i]\in Ker(\alpha)$. Since $A$ is $S$-free with a basis given by $\{(\omega-h_1)^i(\mu-h_2)^j\:|\:0\leq i,j\leq p-1\}$, we get that that $u_i=0$ for $m\leq i\leq p^2$, where $m=p^2-2^{-1}(p-1)p+1$. Let $p\epsilon=\sum_{i=1}^{m-1} v_ix_i$ with the $v_i\in S$ and $x_i$ from the ordered set $T$. Since $p\in S$ is prime we get the following free resolution of $R$ over $S$:
\begin{equation}\label{freeres}
\xymatrix@C+1pc{
 0\ar[r] & S\ar[r]^{\psi^T}
  &
  S^{p^2+1}  \ar[r]^{\alpha}
  & R \ar[r] & 0 
}
\end{equation} 
where $\psi=[v_1 \dots v_{m-1}\;0\dots 0\;-p]$. The above resolution is minimal since $\psi^T(S)\subseteq \mathfrak{m}S^{p^2+1}$, so that $p.d_S(R)=1$. The proof is now complete. 
\QED
\begin{rem}\normalfont \label{pfg}
Let $S$ be an unramified regular local ring of mixed characteristic $p\geq 3$. Note that the ideal $(p,h_1,h_2)\subseteq S$ is a two generated ideal or all of $R$ if and only if the same property holds for the ideal $(p,f,g)\subseteq S$. Similarly, the ideal $(p,h_1,h_2)$ has grade three if and only if the ideal $(p,f,g)$ has the same property.
\end{rem}
\begin{Example}\normalfont\label{ex}
The conditions in \ref{RnotCM} give a non-empty class of non Cohen-Macaulay integral closures $R$. For an example where $Q=(p,h_1,h_2)$ has grade three, consider  $S=\mathbb{Z}[X,Y]_{(3,X,Y)}$ where $X,Y$ are indeterminates over $\mathbb{Z}_{(3)}$. Let \[f=-5X^3+9=X^3+3(3-2X^3)\] \[g=-2Y^3+9=Y^3+3(3-Y^3)\] and $\omega^3=f,\mu^3=g$. Then $f,g$ are square free elements that form a regular sequence in $S$. It is easily checked that $[K:L]=9$ and that this choice satisfies the hypothesis of \ref{RnotCM}, so that $p.d._S(R)=1$.
\par For an example where $Q$ has grade two but $p.d_S(S/Q)=3$, let $S=\mathbb{Z}[X,Y]_{(p,X,Y)}$ for some prime number $p\geq 3$. Set
\[f=(1-p)X^{2p}+p^2=(X^2)^p+p(p-X^{2p})\] \[g=(1+p)(XY)^p+p^2=(XY)^p+p(p+(XY)^p)\]
Then $f,g\in S$ are square free and form a regular sequence in $S$. It is easily verified that $[K:L]=p^2$ and that the choice satisfies the hypothesis of \ref{RnotCM}, so that $p.d._S(R)=1$.
\QED
\end{Example}
\begin{Lemma}\label{T2helplemma}
With established notation, the following holds
\[P^*_A=\langle 1,p^{-1}(\omega-h_1)^{p-1}(\mu-h_2)^{p-1}\rangle_A\]
\end{Lemma}
\demo{Set $P_1:=(p,\omega-h_1)$ and $P_2:=(p,\mu-h_2)$, so that $P^*=P_1^*\cap P_2^*$. Let $\tilde{P_1}:=(p,W-h_1)\subseteq S[\mu][W]$. It is the maximal minors of 
\[   E = \begin{bmatrix} 
    W-h_1 \\
    p\\

    \end{bmatrix}
    \]
    We have $F(W)\in \tilde{P_1}$, $F(W)=a\cdot (-p)+(W^{p-1}+h_1W^{p-2}+\dots+h_1^{p-1})(W-h_1)$. Adjoining the appropriate column of coefficients
     \[   E' = \begin{bmatrix} 
    W-h_1  & a\\
    p& W^{p-1}+\dots+h_1^{p-1}\\
\end{bmatrix}
    \]
From \cite{KU}[Lemma 2.5] $P_1^*=\langle E'_{11}/\delta_1,E'_{22}/\delta_2\rangle_A$ where $E'_{ii}$ and $\delta_i$ denote the image in $A$ of the $(i,i)$-th cofactor of $E'$ and the $i$-th (signed) minor of $E$. Therefore \[P_1^*=\langle 1,p^{-1}(\omega^{p-1}+\dots+h_1^{p-1})\rangle_A\] Identically \[P_2^*=\langle 1,p^{-1}(\mu^{p-1}+\dots+h_2^{p-1})\rangle_A\]
Now consider $y\in P^*=P_1^*\cap P_2^*$. Write for some $\alpha_1,\alpha_2,\beta_1,\beta_2\in A$ \[py=p\alpha_1+\beta_1(\omega^{p-1}+\dots+h_1^{p-1})=p\alpha_2+\beta_2(\mu^{p-1}+\dots+h_2^{p-1})\]
Lifting to $B:=S[W,U]$ and denoting lifts by $\sim$
\[p(\tilde{\alpha_1}-\tilde{\alpha_2})+\tilde{\beta_1}(W^{p-1}+\dots+h_1^{p-1})-\tilde{\beta_2}(U^{p-1}+\dots+h_2^{p-1})\in (F(W),G(U))\]
Writing $W^{p-1}+\dots+h_1^{p-1}=(W-h_1)^{p-1}+p\cdot C_1'$ (respectively for $U^{p-1}+\dots+h_2^{p-1}$), 
\[\tilde{\beta_1}(W-h_1)^{p-1}-\tilde{\beta_2}(U-h_2)^{p-1}\in (p,F(W),G(U))\]
This gives $\tilde{\beta_1}\in (p,W-h_1,(U-h_2)^{p-1})$. Since $1\in P^*$ and $(\omega-h_1)(\omega^{p-1}+\dots+h_1^{p-1})\in pA$, we get $P^*\subseteq \langle 1,p^{-1}(\omega-h_1)^{p-1}(\mu-h_2)^{p-1}\rangle_A$. Since the reverse inclusion is obvious, the proof is complete. 
\QED}

\begin{Theorem}\label{T2}
Let $(S,\mathfrak{m})$ be an unramified regular local ring of mixed characteristic $p\geq 3$. Then
\begin{enumerate}
    \item $R$ is Cohen-Macaulay if
    \begin{enumerate}
        \item At least one of $S[\omega], S[\mu]$ is not integrally closed.
        \item $S[\omega],S[\mu]$ are integrally closed and $fg^i\notin S^{p\wedge p^2}$ for all $1\leq i\leq p-1$.
    \end{enumerate}
    \item Let $S[\omega],S[\mu]$ be integrally closed and $fg\in S^{p\wedge p^2}$. Then $R$ is Cohen-Macaulay if and only if $Q:=(p,f,g)\subseteq S$ is a two generated ideal or all of $S$. Moreover, $p.d_S(R)\leq 1$ and $\nu_S(R)\leq p^2+1$.
    \item If $Q:=(p,f,g)\subseteq S$ has grade three, $R$ admits a birational maximal Cohen-Macaulay module.
\end{enumerate} 

\end{Theorem}
\demo
We have shown 1(a) in $\ref{T1}$ and 1(b) in \ref{T1'}. The proof of (2) follows from \ref{RnotCM}, \ref{lem2Rgen} and \ref{pfg}. 
\par Now assume $Q$ has grade three. From part (1), we may assume that $S[\omega], S[\mu]$ are integrally closed and $fg^i\in S^{p\wedge p^2}$ for some (unique) $1\leq i\leq p-1$. From $\ref{P10}$, $I^*=R$ for $I:=pA+(\omega\mu^i-h_1h_2^i)\cdot P^{p-2}$. Set $M:=(IP)^*$. Then $M$ is an $R$-module since $(A:_K IP)=((A:_K I):_K P)=(R:_K P)$. We will show $depth_S(M)=d$, so that $M$ is an MCM module over $R$. By definition \[M=(IP)^*=(p\cdot P+(\omega\mu^i-h_1h_2^i)\cdot P^{p-1})^*=F_1\cap F_2.\] 
where $F_1=p^{-1}P^*$ and $F_2=(\omega\mu^i-h_1h_2^i)^{-1}(P^{p-1})^*$. This is because for ideals $H,N\subseteq A$, $(A:_K H+N)=(A:_KH)\cap (A:_KN)$ as $A$-modules. Now $A/P\simeq S/pS$ as $S$-modules, therefore by the depth lemma $P$ is $S$-free. By $\ref{Gor}$, $Hom_A(P,A)\simeq Hom_S(P,S)$ as $S$-modules and hence $P^*$ is Cohen-Macaulay. On the other hand, since $(P^{p-1})^*$ and $(P^{(p-1)})^*$ are birational $S_2$ modules that agree in codimension one, we have $(P^{p-1})^*=(P^{(p-1)})^*$. From \ref{lem1}(2) and \ref{P9} we then have that $(P^{p-1})^*$ is Cohen-Macaulay. Therefore $F_1$ and $F_2$ are Cohen-Macaulay since $F_1\simeq P^*$ and $F_2\simeq (P^{p-1})^*$ as $A$-modules and $S$-modules. We have the natural short exact sequence of $S$-modules
\begin{equation}\label{depth1}
\xymatrix@C+1pc{
 0\ar[r] & F_1\cap F_2 \ar[r]
  &
  F_1\oplus F_2  \ar[r]
  & F_1+F_2 \ar[r] & 0 
}
\end{equation}    
To complete the proof it suffices to show that $depth_S(F_1+F_2)\geq d-1$. Set \[\mathscr{F}:=p(\omega\mu^i-h_1h_2^i)\cdot (F_1+F_2)= \mathscr{F}_1+\mathscr{F}_2.\] where $\mathscr{F}_1:=(\omega\mu^i-h_1h_2^i)P^*$ and $\mathscr{F}_2:=p(P^{p-1})^*$. Clearly $F_1+F_2\simeq \mathscr{F}$ as $A$-modules and hence as $S$-modules. From \ref{lem1}(1), $\mathscr{F}_2=(p)+P^{p}$ and from \ref{T2helplemma}
\begin{equation}\label{F1A-mod}
\mathscr{F}_1=(\omega\mu^i-h_1h_2^i,p^{-1}(\omega\mu^i-h_1h_2^i)(\omega-h_1)^{p-1}(\mu-h_2)^{p-1})A.\end{equation}
Set $m:=\omega\mu^i-h_1h_2^i$. We make the following two claims:
\begin{enumerate}
    \item $\mathscr{F}=\mathscr{F}_2+(m)$.
    \item $(\mathscr{F}_2:_Am)=(p)+P^{p-1}$.
\end{enumerate}
Assume both claims hold. Since $(\mathscr{F}_2:_Am)\simeq \mathscr{F}_2\cap (m)$ as $A$-modules and hence $S$-modules, we have a natural short exact sequence of $S$-modules
\begin{equation}\label{depth2}
\xymatrix@C+1pc{
 0\ar[r] & (p)+P^{p-1} \ar[r]
  &
  \mathscr{F}_2\oplus (m)  \ar[r]
  & \mathscr{F} \ar[r] & 0 
}
\end{equation}  
If $depth_S((p)+P^{p-1})=d$, then $depth_S(\mathscr{F})\geq d-1$ and we are done. But $depth_S((p)+P^{p-1})=d$ if and only if $A/((p)+P^{p-1})$ is Cohen-Macaulay. For $B:=S[W,U]_{(\mathfrak{m},W-h_1,U-h_2)}$ we have as $B$-modules \[A/((p)+P^{p-1})\simeq B/((p)+(W-h_1,U-h_2)^{p-1}).\]  Since $B/pB$ is regular local and any power of a complete intersection $B$-ideal is perfect, we are through. Therefore only the claims remain to be proved.
\par
Set $\mathscr{Q}:=\mathscr{F}_2+(m)$.
For claim (1), from (\ref{F1A-mod}) we only need to show 
\[s:=p^{-1}m(\omega-h_1)^{p-1}(\mu-h_2)^{p-1}\in \mathscr{Q}.\]
Since $fg^i\in S^{p\wedge p^2}$, we get $ah_2^p+ibh_1^p\in pS$. Moreover, $grade(Q)=3$ implies $a-qh_1^p\in pS$ and $b+i^{-1}qh_2^p\in pS$ for some $q\in S$.
Write 
\begin{equation}\label{identity1}
m=(\omega-h_1)(\mu^i-h_2^i)+h_2^i(\omega-h_1)+h_1(\mu^i-h_2^i).
\end{equation}
and recall that $(\omega-h_1)^p=p(a-c_1'(\omega-h_1))$ and $(\mu-h_2)^p=p(b-c_2'(\mu-h_2))$. Then
 \begin{align}\label{claim1proof}
 \begin{split}
 s  &\equiv p^{-1}[h_2^i(\omega-h_1)+h_1(\mu^i-h_2^i)](\omega-h_1)^{p-1}(\mu-h_2)^{p-1}\; \text{mod}\: \mathscr{Q}
 \\
 &\equiv ah_2^i(\mu-h_2)^{p-1}+bh_1(\mu^{i-1}+\dots+h_2^{i-1})(\omega-h_1)^{p-1}\; \text{mod}\: \mathscr{Q}
 \\
 &\equiv ah_2^i(\mu-h_2)^{p-1}+ibh_1h_2^{i-1}(\omega-h_1)^{p-1}\; \text{mod}\: \mathscr{Q}
 \\
 &\equiv qh_1h_2^i[h_1^{p-1}(\mu-h_2)^{p-1}-h_2^{p-1}(\omega-h_1)^{p-1}]\; \text{mod}\: \mathscr{Q}.
 \end{split}
 \end{align}
Now $(\omega\mu^i-h_1h_2^i)\cdot P^{p-2}\subseteq \mathscr{Q}$, (\ref{identity1}) and $P^p\subseteq \mathscr{Q}$ imply 
\begin{equation*}\label{colon1}
h_2^i(\omega-h_1)+h_1(\mu^i-h_2^i)\in (\mathscr{Q}:_A P^{p-2}).
\end{equation*}
Therefore for all $0\leq j\leq p-2$
\begin{equation*}\label{idealmembership7}
   \begin{aligned}
   &(h_2^i(\omega-h_1)+h_1(\mu^i-h_2^i))\cdot h_2^{p-2-j}h_1^j(\omega-h_1)^{p-2-j}(\mu-h_2)^j\\
   =\: &h_2^{i-1}[(h_2(\omega-h_1))^{p-j-1}(h_1(\mu-h_2))^j]+
   \\
   &(h_1(\mu-h_2))^{j+1}(\mu^{i-1}+\dots +h_2^{i-1})(h_2(\omega-h_1))^{p-2-j}\in \mathscr{Q}.
   \end{aligned}
\end{equation*}
Thus
\begin{equation*}
h_2^{i-1}[(h_2(\omega-h_1))^{p-j-1}(h_1(\mu-h_2))^j]\equiv -ih_2^{i-1}(h_1(\mu-h_2))^{j+1}(h_2(\omega-h_1))^{p-2-j}\; \text{mod}\: \mathscr{Q}.
\end{equation*}
It then follows that for any $1\leq k\leq p-j-1$
\begin{equation*}
h_2^{i-1}[(h_2(\omega-h_1))^{p-j-1}(h_1(\mu-h_2))^j]\equiv (-i)^kh_2^{i-1}(h_1(\mu-h_2))^{j+k}(h_2(\omega-h_1))^{p-j-1-k}\; \text{mod}\: \mathscr{Q}.
\end{equation*}
In particular for $j=0$ and $k=p-1$ we get
\begin{equation}\label{finalmembership}
h_2^{i-1}(h_2(\omega-h_1))^{p-1}\equiv h_2^{i-1}(h_1(\mu-h_2))^{p-1}\; \text{mod}\: \mathscr{Q}.
\end{equation}
Combining (\ref{finalmembership}) and (\ref{claim1proof}), we see that $s\in \mathscr{Q}$ and thus claim (1) holds.
\\
To show one containment in claim (2), note that $(p)+P^{p-1}\subseteq (\mathscr{F}_2:_A m)$ since $p+P^p= \mathscr{F}_2$. For the reverse inclusion, consider $y\in (\mathscr{F}_2:_A m)$. Lifting to $B$ and denoting lifts by $\sim$
\begin{equation}\label{claim2-1}
\tilde{y}(WU^i-h_1h_2^i)\in (p,F(W),G(U))+(W-h_1,U-h_2)^p=(p)+(W-h_1,U-h_2)^p.
\end{equation}
Using (\ref{identity1}) we have $\tilde{y}\in (p,(W-h_1)^{p-1},U-h_2)$. Similarly $\tilde{y}\in (p,W-h_1,(U-h_2)^{p-1})$ and from (\ref{regsequence})
\begin{equation}\label{induction1}
    \tilde{y}\in (p,(W-h_1)^{p-1},(U-h_2)^{p-1},(W-h_1)(U-h_2)).
\end{equation}
Now assume for some $1\leq i\leq 2^{-1}(p-1)-1$,
\begin{equation}\label{claim2-2}
\tilde{y}\in (p,(W-h_1)^i(U-h_2)^i)+(W-h_1,U-h_2)^{p-1}.
\end{equation}
Write $\tilde{y}=\alpha\cdot (W-h_1)^i(U-h_2)^i+\beta$ for some $\alpha\in B$ and $\beta\in (p)+(W-h_1,U-h_2)^{p-1}$. From $(\ref{claim2-1})$:
\begin{equation}\label{claim2-3}
 \alpha\cdot(W-h_1)^i(U-h_2)^i(WU^i-h_1h_2^i)\in (p)+(W-h_1,U-h_2)^p.
\end{equation}
Using the regular sequence $(p,(U-h_2)^{i+1},(W-h_1)^{i+1},h_2^i)\subseteq B$, we get \[\alpha\in (p,(W-h_1)^{p-2i-1},U-h_2).\] 
Similarly, using the regular sequence $(p,(W-h_1)^{i+1},(U-h_2)^{i+1},h_1(U^{i-1}+\dots +h_2^{i-1}))\subseteq B$ we get $\alpha\in (p,W-h_1,(U-h_2)^{p-2i-1})$. Thus by (\ref{regsequence}):
\[\alpha\in(p,(W-h_1)^{p-2i-1},(U-h_2)^{p-2i-1},(W-h_1)(U-h_2))\]
and hence \[\tilde{y}\in (p,(W-h_1)^{i+1}(U-h_2)^{i+1})+(W-h_1,U-h_2)^{p-1}.\]
Thus starting from (\ref{induction1}) we may induct on $i$ to get 
\[\tilde{y}\in (p,(W-h_1)^{2^{-1}(p-1)}(U-h_2)^{2^{-1}(p-1)})+(W-h_1,U-h_2)^{p-1}=(p)+(W-h_1,U-h_2)^{p-1}.\]
This shows $(\mathscr{F}_2:_A m)=(p)+P^{p-1}$ and all claims have been proved. Thus $R$ admits a birational MCM module.

\QED 
\begin{rem}\normalfont
If $grade(Q)=2$ and $p.d_S(S/Q)=3$ in the context of \ref{T2}(3), we are not able to construct a birational MCM module over $R$ at present. However, if we allow an extension of the quotient field, then constructing an MCM module over $R$ may be possible in this case, see \cite{DKPS}.
\end{rem}
\begin{rem}\normalfont\label{VB}
By a vector bundle on the punctured spectrum of a regular local ring $(S,\mathfrak{m})$ or simply a bundle on $S$ we mean a finitely generated reflexive $S$-module $M$ such that $M_P$ is $S_P$-free for all non maximal ideals $P\subseteq S$. One could use \ref{T2}(2) to generate examples of non-trivial bundles $M$ on localizations of polynomial rings or power series rings over $\mathbb{Z}_{(p)}$ of dimension $d$ at least three such that $rank_S(M)=p^2+d-3$. Moreover, these bundles would satisfy $p.d._S(M)=1$.
\par 
Let $d\geq 3$ and $(T,\mathfrak{n})$ be a $d$-dimensional unramified regular local ring of mixed characteristic $p$. Choose $(S,\mathfrak{m})\subseteq T$ a three dimensional subring of $T$ that is an unramified regular local ring of mixed characteristic $p$ and a quotient of $T$ by a regular sequence (such a choice is possible for example when $T$ is a localization of a polynomial ring over $\mathbb{Z}_{(p)}$). Let $S\subseteq (E,\mathfrak{n}')\subseteq T$ be such that $E$ is regular local and $S=E/(t)$ for some $0\neq t\in E$. Using $\ref{T2}(2)$, with the base ring as $S$, construct $R$ such that it is not $S$-free. Choose a minimal $S$-free resolution 
\begin{equation}\label{freeres2}
\xymatrix@C+1pc{
 0\ar[r] & S\ar[r]^{\psi^T}
  &
  S^{p^2+1}  \ar[r]
  & R \ar[r] & 0 
}
\end{equation}
Let $M'$ be the cokernel of the $E$-matrix $\phi:=\left[
\begin{array}{ccc}
\psi^T & \vline & t
\end{array}
\right]$
\begin{equation}\label{freeres3}
\xymatrix@C+1pc{
 0\ar[r] & E\ar[r]^{\phi}
  &
  E^{p^2+2}  \ar[r]
  & M' \ar[r] & 0 
}
\end{equation}
so that $p.d._E(M)=1$. The ideal of maximal minors of $\psi^T$ is $\mathfrak{m}$-primary since $R$ is a bundle over $S$. Therefore the ideal of maximal minors of $\phi$ is $\mathfrak{n}'$-primary and hence it is free on the punctured spectrum of $E$. Proceeding this way, we can construct a finite module $M$ over $T$ that is free on the punctured spectrum of $T$ and $p.d._T(M)=1$. Moreover, since $M$ is an $S_2$ $T$-module, it is $T$-reflexive (see \cite{bruns_herzog_1998}[Proposition 1.4.1] for example).
\end{rem}

\section*{Acknowledgement}
I would like to thank my Ph.D. advisor Prof. Daniel Katz for suggesting this problem and for his support through the course of this work. I would also like to thank the referee for a careful reading of the paper and suggestions for improvement.
\bibliographystyle{alpha}

\end{document}